\documentclass{gtart_h}

\def\ifplaintex{\expandafter\ifx\csname documentclass\endcsname\relax}

\ifplaintex 
\else
\fi

\def\gt{{\mathsurround=0pt\it $\cal G\mskip-2mu$eometry \&\ 
$\cal T\!\!$opology}}        

\def\gtm{{\mathsurround=0pt\it $\cal G\mskip-2mu$eometry \&\ 
$\cal T\!\!$opology $\cal M\mskip-1mu$onographs}}    

\def\gtp{{\mathsurround=0pt\it $\cal G\mskip-2mu$eometry \&\ 
$\cal T\!\!$opology $\cal P\!$ublications}}  

\def\recd{{\small Received:\qua\receiveddate\ifx\reviseddate\relax
\else\qquad Revised:\qua\reviseddate\fi\par}} 

\def\volumenumber#1{\def\thevolumenumber{#1}}
\def\volumeyear#1{\def\thevolumeyear{#1}}
\def\volumename#1{\def\thevolumename{#1}}
\def\papernumber#1{\def\thepapernumber{#1}}
\def\pagenumbers#1#2{\def\startpage{#1}\def\finishpage{#2}}
\def\published#1{\def\publishdate{#1}}
\def\received#1{\def\receiveddate{#1}}
\def\revised#1{\def\reviseddate{#1}}
\def\accepted#1{\def\accepteddate{#1}}

\def\asciiemail#1{\def\theasciiemail{#1}}

\let\\\par
\let\thevolumenumber\relax\let\thepapernumber\relax
\let\thevolumeyear\relax\let\startpage\relax
\let\finishpage\relax\let\publishdate\relax\let\receiveddate\relax
\let\reviseddate\relax\let\accepteddate\relax\let\theasciititle\relax
\let\theasciiauthors\relax
\let\theasciiabstract\relax

\let\theerratum\relax\let\theasciiemail\relax
\let\theshortauthors\relax\let\theshorttitle\relax

\def\startpage{1}\def\finishpage{15}\def\thepapernumber{77}

\volumenumber{2}
\volumename{Proceedings of the Kirbyfest}
\volumeyear{1999}

\long\def\maketitlep{   

\count0=\startpage

\gtm\nl        
{\small Volume \thevolumenumber: \thevolumename\nl 
\ifx\theerratum\relax\else Erratum \erratumnumber\nl\fi
Pages \startpage--\finishpage\nl}

\vglue 0.1truein   

{\parskip=0pt\leftskip 0pt plus 1fil\def\\{\par\smallskip}{\ifplaintex\large
\else\Large\fi\bf\thetitle}\par\medskip}   
\vglue 0.05truein 

{\parskip=0pt\leftskip 0pt plus 1fil\def\\{\par}{\sc\theauthors}
\par\medskip}%
 
\vglue 0.03truein 

{\small\leftskip 25pt\rightskip 25pt{\bf Abstract}\stdspace\theabstract

{\bf AMS Classification}\stdspace\theprimaryclass
\ifx\thesecondaryclass\relax\else; \thesecondaryclass\fi\par
{\bf Keywords}\stdspace \thekeywords\par}\vglue 7pt

}   

\font\phead=cmsl9 scaled 950
\font=cmsl9 scaled 1050
\font\pnum=cmbx10 scaled 913
\font\lnum=cmbx10 
\font\pfoot=cmsl9 scaled 950
\font=cmsl9 scaled 1050
\ifplaintex
\headline{\vbox to 0pt{\vskip -4.5mm\line{\small\phead\ifnum
\count0=\startpage ISSN 1464-8997 (on line)
1464-8989 (printed) \hfill {\pnum\folio}\else\ifodd\count0\def\\{ }%
\ifx\theshorttitle\relax\thetitle\else\theshorttitle\fi\hfill{\pnum\folio}
\else\def\\{ and }{\pnum\folio}\hfill\ifx\theshortauthors\relax\theauthors
\else\theshortauthors\fi\fi\fi}\vss}}
\footline{\vbox to 0pt{\vglue 0mm\line{\small\pfoot\ifnum\count0=\startpage
Published \publishdate:\qua\copyright\ \gtp\hfill\else
\gtm, Volume \thevolumenumber\ (\thevolumeyear)\hfill\fi}\vss
}}
\else
\makeatletter
\def\@oddhead{{\small\lhead\ifnum\count0=\startpage ISSN 1464-8997 (on line)
1464-8989 (printed) \hfill {\lnum\number\count0}\else\ifodd\count0
\def\\{ }\ifx\theshorttitle\relax \thetitle \else\theshorttitle\fi\hfill
{\lnum\number\count0}\else\def\\{ and }{\lnum\number\count0}
\hfill\ifx\theshortauthors\relax 
\theauthors\else\theshortauthors\fi\fi\fi}}\def\@evenhead{@oddhead}
\def\@oddfoot{\small\lfoot\ifnum\count0=\startpage Published \publishdate:\qua\copyright\ \gtp\hfill\else
\gtm, Volume \thevolumenumber\ (\thevolumeyear)\hfill\fi}
\def\@evenfoot{@oddfoot}
\makeatother
\fi

\let\maketitlepage\maketitlep
\let\makeshorttitle\maketitlepage
\let\maketitle\maketitlepage

\newwrite\gtoutfile
\long\gdef\makeheadfile{  
{\def\\{, }\def\s{ }
\immediate\openout\gtoutfile head.xxx
\immediate\write\gtoutfile{Proxy-for: \ifx\theasciiauthors\relax
\theauthors\else\theasciiauthors\fi\s<\ifx\theasciiemail\relax\theemail\else\theasciiemail\fi>}
\immediate\write\gtoutfile{\noexpand\\}
\immediate\write\gtoutfile{Authors: \ifx\theasciiauthors\relax
\theauthors\else\theasciiauthors\fi}
{\def\\{ }\immediate\write\gtoutfile{Title: \ifx\theasciititle\relax
\thetitle\else\theasciititle\fi}}
\immediate\write\gtoutfile{Subj-class: GT or SG, GR etc}
\immediate\write\gtoutfile{MSC-class: \theprimaryclass\ifx\thesecondaryclass\relax\else, \thesecondaryclass\fi}
\immediate\write\gtoutfile{Journal-ref: Geom. Topol. Monogr. \thevolumenumber\s
(\thevolumeyear) \startpage-\finishpage}
\immediate\write\gtoutfile{Comments: Published by Geometry and Topology Monographs at}
\immediate\write\gtoutfile{\s\s\s  http://www.maths.warwick.ac.uk/gt/GTMon\thevolumenumber/paper\thepapernumber.abs.html}
\immediate\write\gtoutfile{\noexpand\\}
\immediate\write\gtoutfile{}
\ifx\theasciiabstract\relax
\immediate\write\gtoutfile{\theabstract}\else
\immediate\write\gtoutfile{\theasciiabstract}\fi
\immediate\write\gtoutfile{}
\immediate\write\gtoutfile{\noexpand\\}
\immediate\write\gtoutfile{}
\immediate\closeout\gtoutfile}}  

\def\maketitlepage{\maketitlep\makeheadfile}
\let\makeshorttitle\maketitlepage
\let\maketitle\maketitlepage

\volumenumber{7}
\volumename{Proceedings of the Casson Fest}
\volumeyear{2004}
\papernumber{9}
\pagenumbers{213}{233}
\received{15 October 2003}
\revised{2 January 2004}
\accepted{10 December 2003}
\published{19 September 2004}

\usepackage{amsmath,amssymb,graphicx,verbatim}

\newtheorem{theorem}{Theorem}
\newtheorem{corollary}[theorem]{Corollary}

\newcommand{\remark}[1]{\bigskip\noindent
{\bf Remark}\qua #1 \bigbreak}
\newcommand{\remarks}[1]{\bigskip\noindent
{\bf Remarks}\qua #1 \bigbreak}
\newcommand{\examples}[1]{\bigskip\noindent
{\bf Examples}\qua #1 \bigbreak}
\newcommand{\bz}{\ensuremath{\mathbb Z}}
\newcommand{\bn}{\ensuremath{\mathbb N}}
\newcommand{\p}{\ensuremath{\mathcal{P}}}

\newcommand{\n}{\ensuremath{\mathcal{NP}}}

\newcommand{\sign}{\textup{sign}}
\newcommand{\lk}{\textup{lk}}
\newcommand{\rank}{\textup{rk}}
\newcommand{\smooth}{\place{0}{-2}{$\frown$}\place{1}{3}{$\smile$}\quad\ }
\newcommand{\foot}{\setcounter{footnote}{1}\footnote}

\newcommand{\place}[3]
  {\text{\kern#1pt
     \smash{\raise#2pt\hbox{\rlap{#3}}}
   \kern-#1pt\kern-1ex}}
\newcommand{\fig}[2]{\includegraphics[scale=#2]{#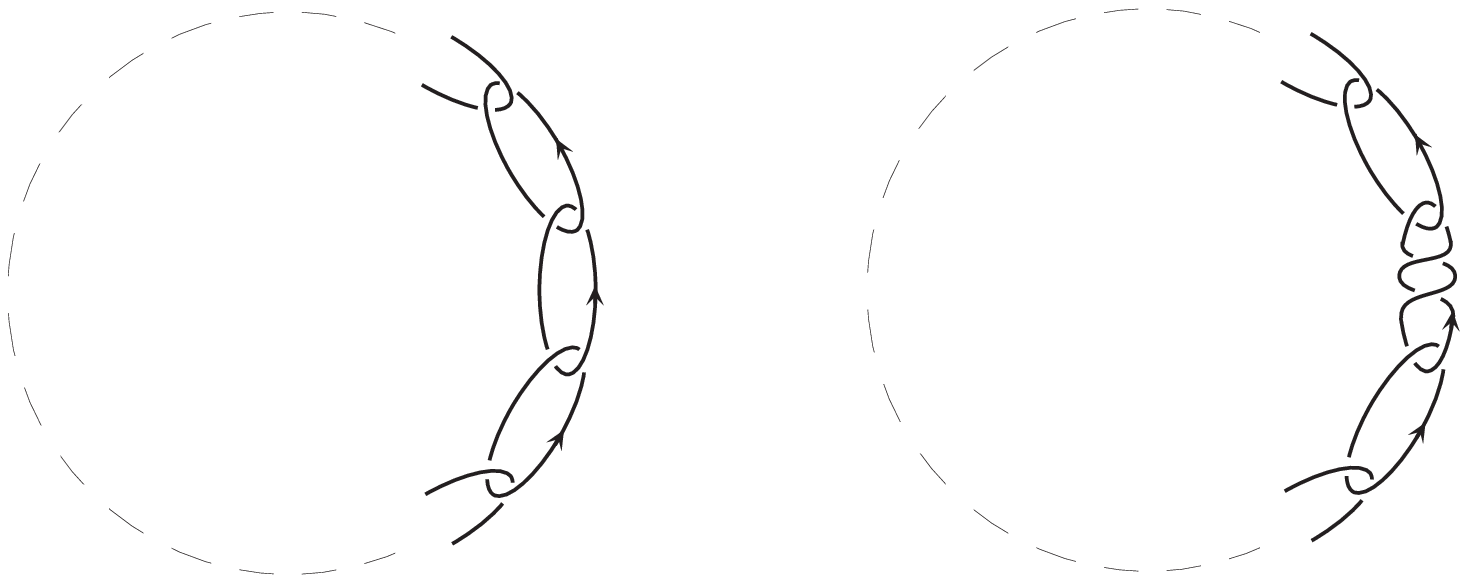}}
\newcommand{\figno}[1]
  {\medskip\centerline{\small Figure #1}\medskip}

\numberwithin{equation}{section}
\numberwithin{theorem}{section}

\begin{document}

\title{Local surgery formulas for quantum invariants\\and the Arf invariant}

\authors{Robion Kirby\\Paul Melvin}

\address{University of California, Berkeley, CA 94720, USA\\Bryn Mawr 
College, Bryn Mawr, PA 19010, USA}

\gtemail{\mailto{kirby@math.berkeley.edu}, \mailto{pmelvin@brynmawr.edu}}
\asciiemail{kirby@math.berkeley.edu, pmelvin@brynmawr.edu}

\begin{abstract}
A formula for the Arf invariant of a link is given in terms of the
singularities of an immersed surface bounded by the link.  This is
applied to study the computational complexity of quantum invariants of
3--manifolds.
\end{abstract}


\maketitlepage

\setcounter{section}{-1}

\section{Introduction}

The quantum 3--manifold invariant of Witten \cite{witten} and
Reshetikhin--Turaev \cite{rt} with gauge group $SU(2)$ at the fourth
root of unity is given by the formula \cite{kme}
$$\tau_4 (M) = \sum_{\theta} \omega^{\mu(M_\theta)}$$
where $\omega$ is a primitive sixteenth root of unity, and the sum is over
all spin structures $\theta$ on the closed oriented 3--manifold $M$.
Here $\mu(M_{\theta})$ is Rokhlin's invariant of $M$ with its spin
structure $\theta$, that is, the signature modulo 16 of any compact spin
4--manifold with spin boundary $M_{\theta}$.  The set of spin structures on
$M$ is parametrized by $H^1(M; \bz_2)$, so at first sight the complexity of
computing $\tau_4 (M)$ grows {\sl exponentially} with $b_1(M) = \rank\,H^1(M;
\bz_2)$.

This note originated when Mike Freedman, motivated by the \p\ versus
\n\ problem in theoretical computer science, observed that the 
formulas in our paper \cite{kme} lead to a {\sl polynomial time} 
algorithm for the computation of $\tau_3$, and asked us what 
difficulties arise in trying to find such an algorithm to evaluate 
$\tau_4$.
As it turns out the computation of $\tau_4$ is \n--hard (and 
conjecturally not even in \n) as we shall explain in section 2, although a 
polynomial time algorithm exists for the restricted class of 
3--manifolds of ``Milnor degree'' greater than three.

We thank Freedman for several discussions on this topic which
led naturally to the ``local'' formulas given below.  We are also
grateful to L\'aszl\'o Lov\'asz for sharing his computational complexity
insights with us.

In the process of investigating this complexity question, we found a new
formula for the Arf invariant of a classical link $L$ in terms of data derived
from an immersed surface $F$ bounded by $L$ whose singularities $S$ 
are {\sl internal}, ie away from $\partial F$.  This formula, 
discussed in section 4 after some algebraic preliminaries in section 3, depends 
only on linking numbers of curves near $\partial F\cup S$, and on the 
Arf invariants (or Brown invariants if $F$ is nonorientable) of 
quadratic forms defined on $H_1(F)$.  For example, if $F$ is a union 
of Seifert surfaces $F_i$ for the individual components $L_i$ of $L$, 
then the formula is expressed in terms of the Arf invariants of the 
$L_i$, the linking numbers between the $L_i$, $F_i\cap F_j$ and their 
push-offs, and the total number of triple points $\cup (F_i \cap F_j 
\cap F_k)$.

The authors thank the National Science Foundation and the Microsoft
Research Group for support.

\section{Local surgery formulas}

It was observed by Casson (see \cite{fk}) and independently by Rokhlin
\cite{rokhlin} that $\mu(M_\theta)$ can be computed using {\sl any} compact
oriented 4--manifold $W$ bounded by $M$ by
$$
\mu (M_{\theta}) = \sigma (W) - F\cdot F + 8\alpha(F) \pmod{16}.
$$
Here $F\subset W$ is an oriented {\sl characteristic surface} for 
$\theta$ --- meaning $\theta$ extends over $W\setminus F$ but not 
across any component of $F$ --- with self intersection $F\cdot F$, and 
$\alpha(F) \in \bz_2$  is the Arf
invariant of a suitable quadratic form on $H_1(F;\bz_2)$.   (See the
appendices of \cite{rs} or \cite{kir} for generalities on the Arf
invariant.)  If $F$ is nonorientable, there is an analogous formula due to
Guillou and Marin \cite{gm}, replacing $8\alpha(F)$ by $2\beta(F)$ where
$\beta(F)\in \bz_8$ is Ed Brown's generalization of the Arf invariant
\cite{brown1}.

In particular, $M$ can be described as the boundary of a 4--manifold 
$W$ obtained from the 4--ball by adding 2--handles along a framed link 
$L$ in $S^3$.   Then the spin structures on $M$ correspond to {\it 
characteristic sublinks} $C$ of $L$, that is sublinks $C$ satisfying
$C\cdot L_i \equiv L_i \cdot L_i \pmod2$
for all components $L_i$ of $L$.  Here \,$\cdot$\, denotes linking or
self-linking number, ie framing (see for example \cite[Appendix
C]{kme}).  Note that linking numbers are only defined for {\sl oriented}
links, so we fix an orientation on $L$; the family of characteristic
sublinks of $L$ is independent of this choice, but we shall need this 
orientation for
other purposes below.

For any given characteristic sublink $C$ of $L$, an associated
characteristic surface $F\subset W$ can be constructed by taking the union
of an oriented Seifert surface for $C$, with its interior pushed into
$B^4$, with the cores of the 2--handles attached to $C$.  The choice of
Seifert surface is immaterial since the invariants $F\cdot F$ and
$\alpha(F)$ that appear in the formula for the $\mu$-invariant are
independent of this choice; indeed, $F\cdot F = C\cdot C$, the sum of all
the entries in the linking matrix of $C$, and $\alpha(F) = \alpha(C)$, the
Arf invariant of the oriented proper link $C$.  (Recall that a link is {\sl
proper} if each component evenly links the union of the other
components.)  It follows that
$$\tau_4(M) = \omega^{\sigma(L)}\sum_C (-1)^{\alpha(C)}\omega^{-C\!\cdot C}$$
where $\sigma(L)$ is the signature of the linking matrix of $L$, and the sum
is over all characteristic sublinks $C\subset L$.  Since there are
$2^{b_1(M)}$ such sublinks, this yields an exponential time evaluation of
$\tau_4$.

In fact the exponential nature of this formula is due solely to the Arf
invariant factors, for without these, the formula could be evaluated in
polynomial time.  To see this note that the linking matrix of $L$
can be stably diagonalized over $\bz$ (eg by \cite{oht}), which
corresponds to adjoining a suitably framed unlink to $L$ and then sliding
handles \cite{kir}.  Once $L$ has been diagonalized, its characteristic
sublinks are exactly those that include all the odd-framed components.
It follows that if there are $b_i$ components of $L$ with framings
congruent modulo 16 to $2i$, then
$$\sum_C\omega^{-C\cdot C} = \omega^{-s} \prod_{0\le i\le7} (1 +
  \omega^{-2i})^{b_i}$$
where $s$ is the sum of all the odd framings in $L$.

Unfortunately the Arf invariant of a proper link $C$ is global
in the sense that its value depends simultaneously on {\sl all} the
components of $C$.  For example the circular daisy chains in Figure 1a and
1b have different Arf invariants but identical families of sublinks
(excluding the whole link).  This casts doubt on the existence of a 
polynomial time algorithm
for computing $\tau_4$.  However if $C$ is {\it totally proper} --- meaning all
the pairwise linking numbers of $C$ are even --- then there exist {\it local
formulas} for $\alpha(C)$, ie formulas that depend only on the sublinks
of $C$ with $k$ or fewer components for some $k$ (such a formula is called
{\it $k$--local}).  It is reasonable to attempt to exploit such local
formulas in the search for an optimal algorithm for computing $\tau_4$.

\centerline{\fig{1}{.6}}
\centerline{\small(a) \kern 160pt (b)}
\figno{1}

One such formula arises from the work of Hoste \cite{hos1} and 
Murakami \cite{mur}.  They showed (independently) that the Arf 
invariant of a totally proper link $C$ can be written as a sum
$$\alpha(C) = \sum_{S<C} c_1(S) \pmod2$$
over all sublinks $S$ of $C$.\foot{This was first observed for knots
by Kauffman \cite{kauffman}, following Levine \cite{levine}, and for 
2--component links by Murasugi \cite{murasugi}.}  Here $c_1(S)$ is the 
coefficient of $z^{s+1}$ in the
Conway polynomial of the $s$--component link $S$.  It is known that
$c_1(S)=0$ if $s>3$ (see for example \cite{hos2}) and so this formula is in
fact 3--local.

This formula can be expressed in geometric terms using familiar 
homological interpretations for the mod 2 reductions of the Conway 
coefficients $c_1(S)$.   As noted above,
$$c_1(S) \equiv  \alpha(S) \pmod2$$
if $S$ is a knot.  

If $S$ has two components, then 
$c_1(S)$ is determined by the linking number of the components and an 
unoriented version of the Sato--Levine invariant, as follows.  The 
oriented Sato--Levine invariant \cite{sato} is defined for any 
oriented 2--component {\it diagonal} link (meaning pairwise linking 
numbers vanish) as the self-linking of the curves of intersection of 
any pair of Seifert surfaces for the components that meet 
transversely in their interiors; it was shown equal to $c_1$ by Sturm 
\cite{sturm} and Cochran \cite{cochran}.  This invariant was extended 
to unoriented (totally) proper links $S$ by Saito \cite{saito} by 
allowing nonorientable bounding surfaces for the components, meeting 
transversely in their interiors in a link $C$.  One then defines
$$\lambda(S) = \lk(C,C^\times) \!\!\pmod8 \ \in \ \bz_8$$
where $C^\times$ is the ``quadruple push-off'' of $C$ (the union of the
boundaries of tubular neighborhoods of $C$ in the two surfaces, oriented
compatibly with any chosen orientation on $C$).  This is shown independent
of the choice of bounding surfaces by a standard bordism argument (see
\cite{saito} for details).  It is clearly an even number, and a multiple of
4 for diagonal links; Saito's invariant is $\lambda/2 \in \bz_4$, and the
oriented invariant is $\lk(C,C^\times)/4 \in \bz$.  (In section 3 we
study a closely related invariant $\delta$, which can take on odd values
as well.)  Saito shows that in general $\lambda(S)$ is congruent mod 4
to the linking number $\lk(S)$ of the two components of $S$ (with any
chosen orientation), and that
$$c_1(S) \equiv \textstyle{\frac{1}{4}}(\lambda(S)+\lk(S)) \pmod2.$$
There is also a Sato--Levine invariant for 3--component diagonal links which
counts the number of signed triple points of intersection of three oriented
Seifert surfaces meeting only in their interiors.  This invariant clearly
depends on an orientation and ordering of the components.  In fact it is
equivalent (up to sign) to Milnor's triple linking number \cite{milnor2}
\cite{turaev}, and its {\sl square} is equal to the Conway coefficient $c_1$
for diagonal links \cite{cochran}.  To extend this invariant to totally
proper links, one must reduce mod 2.   Thus we allow nonorientable 
bounding surfaces and then count triple points mod 2.  This gives a 
$\bz_2$--valued
link concordance invariant $\tau$ by the usual bordism argument.  In fact
$\tau$ is a link homotopy invariant (eg by Cochran's argument in
\cite[Lemma 5.4]{cochran}) which coincides with the mod 2 reduction
of Milnor's triple linking number (cf \cite{mm}), and so we shall call
it the {\it Milnor invariant}.  Now it is not hard to show that
$$c_1(S) \equiv \tau(S) \pmod2.$$
Indeed $c_1 \!\!\pmod2$ is a link homotopy invariant for totally
proper links of at least 3 components (by the proof of
\cite[Lemma 4.2]{cochran}), and the congruence above is easily checked
on Milnor's generators for 3--component link homotopy \cite[page 23]{milnor1}
(as in \cite[page 539]{cochran} in the diagonal case).

Putting these geometric evaluations of $c_1$ into the
Hoste--Murakami formula yields the following:

\begin{theorem}\label{arf}
If $C$ is a totally proper link with components $C_1,\cdots,C_n$, then
$$\alpha(C) = \sum_i\alpha(C_i) + \sum_{i<j}
  {\textstyle{\frac{1}{4}}}(\lambda(C_i,C_j) + \lk(C_i,C_j)) +
  \sum_{i<j<k} \tau(C_i, C_j, C_k) \pmod 2$$
where $\alpha$, $\lambda$ and $\tau$ denote the Arf invariant, unoriented
Sato--Levine invariant and Milnor triple point invariant, respectively.
\end{theorem}

In section 4 this theorem will be rederived as a corollary of a more general
result expressing the {\it Brown invariant} of a link in terms of linking
properties of the singularities of any immersed surface that it bounds.

\section{Complexity}

First recall, in rough terms, the complexity 
classes $\p=$ {\it polynomial time} and $\n=$ {\it nondeterministic 
polynomial time} (see \cite{papa} for a more rigorous discussion). 
A computational problem is said to be in $\p$ if it can be solved by 
an algorithm whose run time on any given instance of the problem is 
bounded by a polynomial function of the size of the instance.  If 
{\sl answers} to the problem can be {\sl checked} in polynomial time, 
then it is said to be in $\n$.  Of course any problem in $\p$ is in 
$\n$, but the converse is not known; this is one of the central open 
problems in theoretical computer science.  

There are a number of 
well-known $\n$ problems, such as the travelling salesman and Boolean 
satisfiability (SAT) problems, whose polynomial time solution would 
yield polynomial time solutions for all $\n$ problems, thus showing 
$\p=\n$.  These are called {\it $\n$--complete} problems.  Any problem 
(whether or not in $\n$) whose polynomial time solution would yield 
polynomial time solutions for all $\n$ problems is said to be {\it 
$\n$--hard}.

With the formula in Theorem 1.1 (in fact the diagonal case is all 
that is needed) it is easy to show that the problem of calculating 
$\tau_4(M)$ for all 3--manifolds $M$ is \n--hard.  The idea is to 
construct a class of 3--manifolds indexed by cubic forms over $\bz_2$ 
whose quantum invariants are given by counting the zeros of the 
associated forms, a well-known \n--hard computational problem.  In 
principle this construction goes back to Turaev's realization theorem 
for ``Rokhlin functions" \cite{turaev2}, but it can be accomplished 
more efficiently in the present setting as follows.

As a warmup, start with the 3--manifolds $M_r$ obtained by zero-framed 
surgery on the
links $L_r$ (for $r=1,2,3$) where $L_1$ is the trefoil, $L_2$ is the
Whitehead link, and $L_3$ is the Borromean rings.  Then $M_r$ has $2^r$
spin structures given by the $2^r$ sublinks $C$ of $L_r$.  The 
$\mu$-invariant is zero in all cases except when $C=L_r$, when it is 
8 (coming from the Arf invariant if $r=1$, the
Sato--Levine invariant if $r=2$, and the Milnor invariant if $r=3$). Thus
$\tau_4(M_r) = 2^r-2$.

More generally, for any cubic form
$$c(x_1,\dots,x_n) = \sum_i c_i\,x_i + \sum_{i,j} c_{ij}\,x_ix_j
  + \sum_{i,j,k} c_{ijk}\,x_ix_jx_k$$
in $n$ variables $x_1,\dots,x_n$ over $\bz_2$, let $L_c$ be the framed link
obtained from the zero-framed $n$--component unlink by tying a trefoil 
knot in each component $L_i$ for which $c_i=1$, a Whitehead link into 
any two components $L_i, L_j$ for which $c_{ij}=1$, and Borromean 
rings into any three components $L_i, L_j, L_k$ for which $c_{ijk} = 1$.  
If $M_c$ is the 3--manifold obtained by surgery on 
$L_c$, then the spin structure
corresponding to any of the $2^n$ characteristic sublinks $C \subset L$
has $\mu$--invariant $8 c(x)$, where $x\in\bz_2^{\,n}$ is the
$n$--tuple with 1's exactly in the coordinates corresponding to the
components of $C$.  Thus
$$\tau_4 (M_c)  = \sum_{x\in\bz_2^{\,n}} (-1)^{c(x)} = 2\#c-2^n$$
where $\#c$ denotes the number of zeros of $c$ (ie solutions
to $c(x)=0$).

\begin{theorem}
For any cubic form $c$, the calculation of $\tau_4(M_c)$ is equivalent to
the calculation of the number $\#c$ of zeros of $c$.  The problem
{\rm \#C} of computing $\#c$ for all cubic forms is \n--hard.
\end{theorem}

\begin{proof}
The first statement was proved above, and the last is presumably well
known to complexity theorists.  We thank L.\ Lov\'asz for suggesting the
following argument.

It is a fundamental result in complexity theory that the Boolean
satisfiability decision problem SAT is \n--complete (Cook's Theorem
\cite{cook}) as is its  ``cubic'' specialization 3--SAT.  It follows that
the associated counting problem \#3--SAT is \n--hard.  This problem asks
for the number $\#e$ of solutions to logical expressions in $n$
variables $x_1,\dots,x_n$ of the form
$$e = a_1 \wedge a_2 \wedge \cdots \wedge a_r$$
where each $a_i$ is of the form $(x_j^\pm \vee x_k^\pm \vee x_\ell^\pm)$.
Here, each $x_i$ can take the value T (true) or F (false); $x_i^+ = x_i$
and $x_i^-$ is the negation of $x_i$;  $\vee$ means {\sl or} and $\wedge$
means {\sl and}.  Thus $\#e$ is the number of ways to assign T or F to
each $x_i$ so that the expression $e$ is true.  

To complete the proof of the theorem, it suffices to produce a polynomial
time reduction of the problem \#3--SAT to \#C.  To achieve
this, consider any logical expression $e$ as above, and rewrite it as
a system of cubic equations over $\bz_2$ by setting $T=0$ and $F=1$ and
replacing $x_i^-$ by $1-x_i$.  Thus each $a_i$ becomes an equation, eg
$(x_j^-\wedge x_k^-\wedge x_\ell)$ becomes $(1-x_j)(1-x_k)x_\ell = 0$.
The resulting system of $r$ cubic equations in $n$ variables has exactly
$\#e$ solutions.

Now change this cubic system into a system of $k=2r$ {\sl quadratic}
equations in $m\le n+r$ variables
\begin{equation}
q = \left\{
\begin{array}{l}
q_1(x_1, \ldots x_m) = 0 \\
\vdots \\
q_k(x_1, \ldots x_m) = 0
\end{array}
\right.
\tag{$*$}
\end{equation}
also with exactly $\#e$ solutions.  In particular, replace each cubic
equation by two quadratics, the first assigning a new variable to the
product of any two of the variables in the cubic, and the second obtained
by substituting this into the cubic.  For example $(1-x_j)(1-x_k)x_\ell =
0$ is replaced by $x_{jk} = x_j x_k$ and $(1-x_j- x_k + x_{jk} )x_\ell  =
0$.

Finally, convert $q$ into a cubic equation by introducing $k$
new variables $z_i$:
$$c = \sum_{i=1}^k z_iq_i(x_1, \ldots x_n) = 0.$$
The number of solutions $\#c$ is equal to $2^k\#e + 2^{k-1}(2^m-\#e)
=2^{m+k-1} + 2^{k-1}\#e $, since any solution to $(*)$ allows any of
$2^k$ choices for the $z_i$, and any non-solution to $(*)$ allows only
$2^{k-1}$ choices for the $z_i$.  Thus an algorithm to evaluate $\#c$ would
yield one of the same complexity for $\#e$, and so \#C is at least as hard
as \#3--SAT.
\end{proof}


\begin{corollary} The calculation of $\tau_4(M)$ for all 3--manifolds
is \n--hard.
\end{corollary}

In particular, Theorem 2.1 shows that this 
calculation for the special class of 3--manifolds $M_c$ arising from 
cubic forms $c$ is already \n--hard (and presumably not in \n, cf 
\cite[section 18]{papa}).  

{\it Added in proof}: In fact one 
need only consider the class of 3--manifolds $M_q$ arising from {\sl 
quadratic} forms $q(x_1,\dots,x_n) = \sum_i c_i\,x_i + \sum_{i,j} 
c_{ij}\,x_ix_j$, so the vanishing of the triple linking numbers does not 
reduce the complexity of the calculation if there are still pairwise 
Whitehead linkings.  This follows by essentially the same proof, 
using the surprising result of Valiant \cite{valiant} (brought to our 
attention by Sanjeev Khanna) that \#2--SAT is also \n--hard, although 
2--SAT is in $\p$!

\remark
{There is a 3--manifold invariant that captures the complexity
of the calculation of $\tau_4(M)$, namely the {\it Milnor degree} $d(M)
\in \bn$ introduced in \cite[page 116]{cm}.  This invariant can be defined
by the condition $d(M)>n$ if $M$ can be obtained by surgery on an
integrally framed link whose $\overline\mu$--invariants of order
$\le n$ vanish (where the {\it order} of a $\overline\mu$ invariant is one 
less than its length, eg the order--2 invariants are Milnor's 
triple linking numbers) \cite{milnor2}.  It follows from the 
discussion in section 1 that there is a polynomial time algorithm for 
computing $\tau_4(M)$ for all 3--manifolds of Milnor degree $>3$, and 
from the discussion above that the computation for 3--manifolds of 
Milnor degree $\le3$ is \n--hard.}

\section{The Brown invariant: algebra}

The Brown invariant  \cite{brown2}, which is a generalization of the Arf
invariant, classifies $\bz_4$--enhanced inner product spaces over
$\bz_2$.  There are many excellent treatments of this subject in the
literature (see eg \cite{brown2,pinkall,gm,mat,kt})
but generally in the context of {\sl nonsingular} spaces.  For the
reader's convenience, with apologies to the experts, we give an
exposition which includes the case of singular forms (cf
\cite{kv,kme,gilmer}).

The example to keep in mind is the space $H_1(F)$ with its
intersection pairing, where $F$ is a compact surface with boundary.
(Throughout this paper, $\bz_2$--coefficients will be assumed.)  The
enhancements in this case arise from immersions of $F$ in $S^3$, and
these give rise to Brown invariants of the links on the boundary of 
$F$, as will be discussed in the next section.

\subsection{Enhanced spaces}

Let $V$ be a finite dimensional $\bz_2$--vector space with a possibly
singular inner product $(x,y)\mapsto x\cdot y$.   Then $V$ splits as an
orthogonal direct sum
$$V = U \oplus V^\perp$$
where $\cdot$ is nonsingular on $U$ and vanishes identically on
$V^\perp = \{x\in V \ | \ x\cdot y = 0 \textup{ for all } y\in V\}$.
(For surfaces, the splitting of $H_1(F)$ arises from a decomposition $F =
C\#D$ where $C$ is closed and $D$ is planar, and so $H_1(F)^\perp$ is the
image of the map $H_1(\partial F) \to H_1(F)$ induced by inclusion.)

A standard diagonalization argument shows that $U$ splits as a sum of
indecomposables of one of two types: the 1--dimensional space $P$ defined
by $x\cdot x=1$ on any basis $x$  (corresponding to the real projective
plane) and the 2--dimensional space $T$ defined by $x\cdot x = y\cdot y =
0$ and $x\cdot y=1$ on any basis $x,y$ (the torus).  Similarly $V^\perp$ is
a sum of trivial 1--dimensional spaces $A$ where $x\cdot x = 0$ (a boundary
connected sum of annuli).  Thus $V$ is built from the indecomposables $P$,
$T$ and $A$, given by the matrices
$$P = (1) \qquad
  T = \begin{pmatrix} 0 & 1 \\ 1 & 0 \end{pmatrix} \qquad
  A = (0).$$
The only relations among these spaces follow from the well-known
isomorphism $P\oplus T \cong 3P$ ($= P\oplus P\oplus P$).  Hence $V$ is
uniquely expressible as a sum of copies of $T$ and $A$ if it is {\it even}
(ie $x\cdot x=0$ for all $x\in V$, which corresponds to orientability for
surfaces), and of copies of $P$ and $A$ otherwise.

Now equip $V$ with a $\bz_4$--valued {\it quadratic enhancement}, that
is a function
$$e\co V\to\bz_4$$
satisfying 
$e(x+y) = e(x)+e(y)+2(x\cdot y)$ 
for all $x,y\in V$.

If $e$ vanishes on $V^{\perp}$, then it is called a {\it proper}
enhancement.  The pair $(V,e)$, also denoted $V_e$, is called an {\it
enhanced space}.  Observe that $e$ is determined by its values on a basis
for $V$, and these values can be arbitrary as long as they satisfy $e(x)
\equiv x\cdot x \pmod2$.  Thus there are $2^{\dim V}$ distinct enhancements
of $V$.  However many of these may be isomorphic, and
the Brown invariant
$$\beta\co\{\textup{enhanced spaces}\}\to\bz_8^*,$$
where $\bz_8^* = \bz_8\cup\{\infty\}$, provides a complete isomorphism
invariant.

\subsection{The Brown invariant of an enhanced space}

Let $V_e$ be an enhanced space.  Perhaps the simplest definition of 
the Brown invariant $\beta(V_e)$ is based on the relative values of 
$e_0$ and $e_2$, and of $e_1$ and $e_3$, where $e_i$ denotes the 
number of $x\in V$ with $e(x)=i$, according to
the scheme indicated in Figure 2.\foot{This definition is in the spirit of the
characterization of the Arf invariant for $\bz_2$--valued enhancements $e$
(where $e(x+y)=e(x)+e(y)+x\cdot y$) in terms of the relative values of
$e_0$ and $e_1$ (see eg \cite{browder}):  $\alpha = 0,\ 1$ or
$\infty$ according to whether $e_0-e_1$ is positive, negative or zero.
Observe that such an $e$ can also be viewed as a $\bz_4$--valued enhancement
by identifying $\bz_2$ with $\{0,2\} \subset \bz_4$, and then $\beta =
4\alpha$.}  For example $\beta=7$ iff $e_0>e_2$ and $e_1<e_3$, and
$\beta=\infty$ iff $e_0=e_2$ and $e_1=e_3$.

{\small\vskip .1in
\centerline{$\sign(e_1-e_3)$}
\vskip .05in
\centerline{\fig{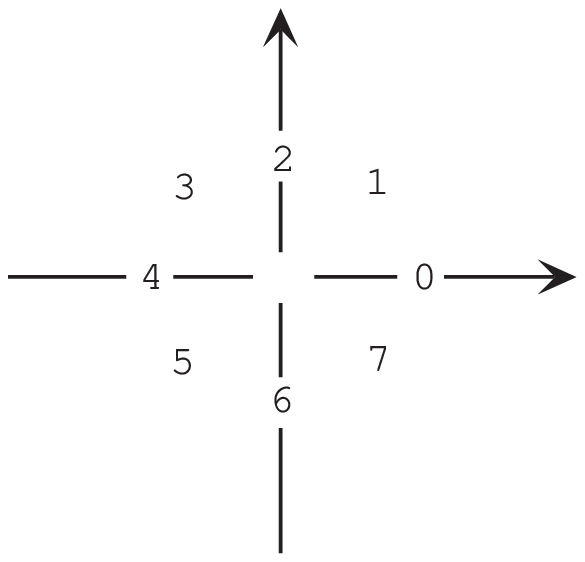}{.8}}
\place{260}{82}{$\sign(e_0-e_2)$}
\place{177}{82}{$\infty$}

\vskip-.2in\figno{2: The Brown invariant}}

Now $T$ has four enhancements $T^{0,0}, T^{0,2}, T^{2,0}, T^{2,2}$ (where
the superscripts give the values on a basis) which fall into two isomorphism
classes, $T_0 = \{T^{0,0}, T^{0,2},T^{2,0}\}$ and $T_4 = \{T^{2,2}\}$ (the
subscripts specify the Brown invariants).  The spaces $P$ and $A$ each have
two nonisomorphic enhancements $P_1,P_{-1}$ and $A_0, A_\infty$ (where once
again the subscripts are the Brown invariants).  Thus $V_e$ decomposes as a
sum of copies of $T_0$, $T_4$, $P_\pm1$, $A_0$ and $A_\infty$, and it is
proper if and only if there are no $A_\infty$ summands.  

The isomorphism $P\oplus T \cong 3P$ above induces isomorphisms 
$$P_{\pm1} \oplus T_0 \cong P_{\pm1} \oplus P_1 \oplus P_{-1}$$
and
$$P_{\pm1} \oplus T_4 \cong 3P_{\mp1},$$
and the 
latter implies the first of the following basic relations (the others 
are left as exercises):
\begin{enumerate}
\item[(a)] $4\,P_1 \cong 4\,P_{-1}$ \,\ and \,\ $2\,T_0 \cong 2\,T_4$
\  (see \cite{brown2} or \cite{mat} for details)
\item[(b)] $P_1 \oplus A_\infty \cong P_{-1} \oplus A_\infty$, \
$T_0 \oplus A_\infty \cong T_4 \oplus A_\infty$ \,\ and \,\ $A_0 \oplus
A_\infty \cong A_\infty \oplus A_\infty$.
\end{enumerate}
It follows from (b) that any two
improper enhancements on $V$ are isomorphic, and from (a) that for
proper enhancements $e$,
\begin{enumerate}
\item[$\bullet$] $V$  even $\implies V_e$ is a sum of copies of
$T_0,T_4,A_0$, with at most one $T_4$
\item[$\bullet$] $V$ odd $\implies V_e$ is a sum of copies of 
$P_{\pm1},A_0$, with at most three $P_{-1}$'s.
\end{enumerate}
In fact these decompositions are unique since the Brown invariant adds under
orthogonal direct sums.  This additivity can be seen using the Gauss sum
$$\gamma(V_e) = \sum_{x\in V} i^{e(x)}$$ which clearly multiplies under
$\oplus$.  One readily computes $\gamma(P_{\pm1}) = 1\pm i$, $\gamma(T_0) =
2$, $\gamma(T_4) = -2$, $\gamma(A_0) = 2$ and $\gamma(A_\infty) = 0$.  It
follows that $\gamma(V_e) = 0$ if $e$ is improper, and by the definition
of $\beta$
$$\gamma(V_e) = \sqrt2^{\,m+n} \, \exp(\pi i\beta(V_e)/4)$$
if $e$ is proper, where $m=\dim V$ and $n=\dim V^{\perp}$.  The additivity
of $\beta$ now follows from the multiplicativity of $\gamma$.

\section{The Brown invariant: topology}

\subsection{The Brown invariant of an immersion}

Let
$$f\co F\looparrowright S^3$$
be an immersion of a compact surface $F$.  The immersion is assumed 
to be {\sl regular}, meaning that the only singularities of $f$ are 
interior transverse double curves with isolated triple points.  Then 
there is an associated quadratic enhancement
$$f_*\co H_1(F)\to\bz_4$$
(recall that $\bz_2$--coefficients are used throughout) which, in rough terms,
counts the number of half-twists modulo 4 in the images of band
neigborhoods of cycles on $F$.  This is defined precisely below.  The {\it
Brown invariant} of $f$ is defined to be the Brown invariant of this 
enhanced space,
$$\beta_f = \beta(H_1(F)_{f_*}).$$
If $f$ is an embedding then it can be identified with its image, and 
we write $\beta(F)$ for $\beta_f$.

To make this precise, we follow an approach suggested by Sullivan
\cite[Example 1.28]{brown2} and later developed by Pinkall 
\cite[section 2]{pinkall} and Siebenmann \cite[Appendix]{gm}.  (Also see 
Guillou and Marin \cite{gm} or Matsumoto \cite{mat} for the analogous 
theory for closed surfaces in simply-connected 4--manifolds.)  Define 
a {\it band} to be a union of annuli
and M\"obius strips, and consider the function
$$\widehat h\co\{\textup{embedded bands in $S^3$}\} \to \bz$$
given by $\widehat h(B) = \lk(C,\partial B)$, where $C$ is the core of the
embedded band $B$ (its zero-section when viewed as an $I$--bundle) and
$\partial B$ is its boundary.  Here $C$ and $\partial B$ should be
oriented compatibly on components, as shown in Figure 3a.

{\small\vskip .1in
\centerline{\fig{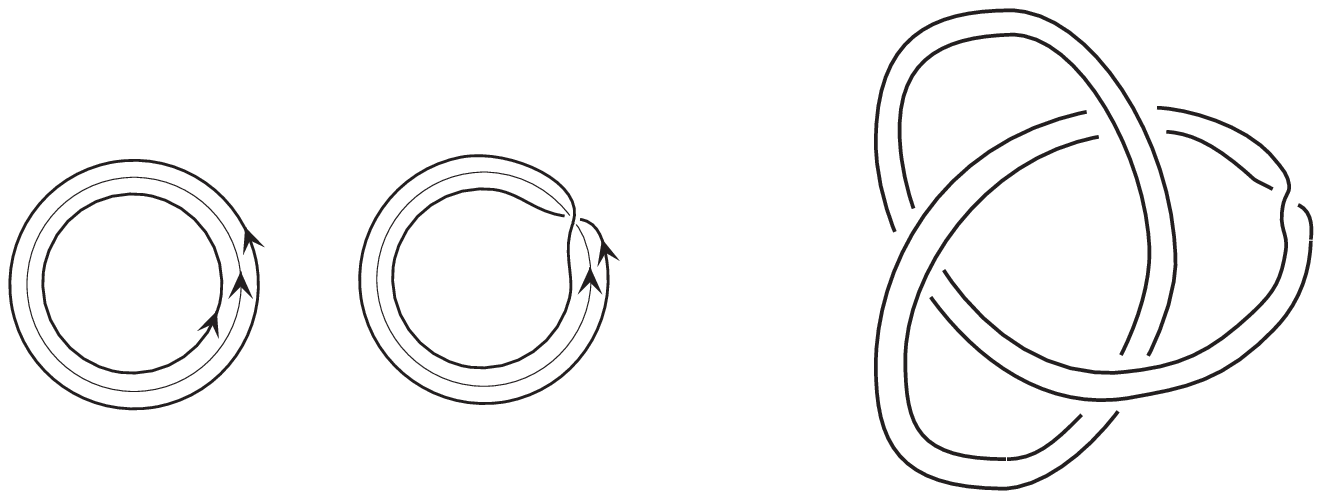}{.6}}
\centerline{(a) compatible orientations \kern 50pt (b) \ $\widehat h =
1+2\cdot3 = 7$}
\figno{3}}

If $B$ is connected, then $\widehat h(B)$ is just the number of
right half-twists in the band relative to the corresponding zero-framed
annular band, computed from a projection as ``twist'' (number of
half-twists) plus twice the ``writhe'' (signed sum of the self-crossings
of the core).  An example is shown in Figure 3b.

Observe that the mod 4 reduction $h(B)$ of $\widehat h(B)$ is unaffected by
linking among the components of $B$, and is in fact invariant under any
regular homotopy of $B$, since a band pass changes $\widehat h$ by 4.  It
follows that there is a well-defined function
$$h\co\{\textup{immersed bands in $S^3$}\} \to \bz_4,$$
which will be called the {\it half-twist map}.  This map is additive 
under unions (meaning $h(B\cup B') = h(B) + h(B')$, where $B$ and $B'$
may intersect) and is a complete regular homotopy invariant for connected
bands.

Now define the enhancement $f_*$ induced by $f$ as follows (being careful,
at least when $f$ is not an embedding, to distinguish between subsets
$S\subset F$ and their images $S' = f(S)\subset S^3$):   For $x\in H_1(F)$,
choose a regularly immersed cycle $C$ in $F$ representing $x$, and set
$$f_*(x) \ = \ h(B') + 2 d(C) \pmod 4$$
where $B$ is an immersed band neighborhood of $C$ (with image $B' \subset
S^3$) and $d(C)$ is the number of double points of $C$ in $F$.

To check that this definition is independent of the choice of $C$, first
observe that small isotopies of $C$ do not change the right hand side.
Thus we may assume that $C$ is transverse to the double curves of $f$,
and that $f$ embeds $B$ onto an immersed band $B'$.  

Now consider the
special case in which $C$ is embedded and null-homologous in $F$, and so
in particular $B'$ is an embedded band neighborhood of $C'$.  We must
show $h(B')=0$.  But $C$ bounds a surface in $F$ whose interior $E$ has
image $E'$ transverse to $C'$ at an {\sl even} number of points (an easy
exercise) and so
$h(B') \equiv \lk(C',\partial B') \equiv 2C'\cdot E' \equiv 0 \pmod4.$

In general, if $C_1$ and $C_2$ are two regular cycles representing $x$,
then after a small isotopy into general position, $C = C_1\cup C_2$ is a
regular null-homologous cycle.  Smoothing crossings $\times
\rightsquigarrow \smooth$ converts $C$ into an embedded cycle
without changing $h+2d$ (each smoothing changes both terms by 2) and so
$h(B') + 2d(C) \equiv 0 \pmod4$ by the special case above.  Since $h$ is
additive, $h(B') = h(B_1') + h(B_2')$, while $d(C) = d(C_1) + d(C_2)
+C_1\cdot C_2$.  Rearranging terms gives
\begin{eqnarray*}
h(B_1')+2d(C_1) &\equiv& (-h(B_2')-2C_1\cdot C_2) - 2d(C_2) \\
&\equiv& h(B_2') + 2d(C_2) \pmod4
\end{eqnarray*}
since $h(B_2') \equiv x\cdot x \equiv C_1\cdot C_2 \pmod2$.  Thus
$f_*$ is well-defined (compare Propositions 1 and 2 in \cite{gm} and
Lemma 5.1 in \cite{mat}).

It is now immediate from the definitions that $f_*$ is quadratic.
Furthermore, it is readily seen that $f_*$ is proper if and only if the
link $L = f(\partial F)$ is proper, ie each component $K$ of $L$ links
$L-K$ evenly.  Indeed, if $K^+$ is a parallel copy of $K$ (the image of a
push off in $F$), then $f_*([K]) = 2\,\lk(K,K^+) = 2 \, \lk(K,L-K)$ (since
$K^+$ and $L-K$ are homologous in $S^3-K$ across $F'$) and so $f_*([K]) =
0$ if and only if $\lk(K,L-K)$ is even.  In this case (when $L$ is proper)
we shall refer to $f$ as a {\it proper immersion}.

\subsection{The Brown invariant of a proper link}

Observe that the Brown invariant of a proper {\sl embedded} surface $F \subset
S^3$ depends only on the {\sl framed} link $L = \partial F$, where the
framing is given by a vector field normal to $L$ in $F$.  (Note that each
component gets an {\sl even} framing since $F$ is proper.)  For if $F'$
is any other surface in $S^3$ bounded by $L$ with the {\sl same} framing,
then the closed surface $S \subset S^4$ obtained from $F\cup F'$ by
pushing int$(F)$ and int$(F')$ to opposite sides of an equatorial $S^3$
has Brown invariant $\beta(S) = \beta(F) - \beta(F')$ and
self-intersection $S\cdot S=0$ (defined to be the twisted Euler class of
the normal bundle of $S$ in $S^4$ when $S$ is nonorientable, cf
\cite{gm}.)  But $\beta(S)=0$ since
$2\beta(S) \equiv \sigma(S^4) - S\cdot S \pmod{16}$
by a theorem of Guillou and Marin \cite{gm}, where $\sigma$ denotes the
signature, and so $\beta(F) = \beta(F')$.

Thus one is led to define the Brown invariant for any {\sl even framed
proper link} $L$  by
$
\beta(L) = \beta(F)
$
where $F\subset S^3$ is any embedded surface bounded by $L$ which induces
the prescribed framing on $L$.  Such a surface can be constructed from an
arbitrary surface bounded by $L$ by stabilizing (adding small half-twisted
bands along the boundary to adjust the framings), and any two have the
same Brown invariant by the discussion above.

If no framing is specified on $L$, then the zero framing is presumed.  In
other words, the {\it Brown invariant of a  proper link} (unframed) is
defined to be the Brown invariant of the link with the zero framing on
each component.  It can be computed from {\sl any} embedded surface $F$
bounding $L$ (possibly nonorientable) by the formula
$$
\beta(L) = \beta(F) - \phi(F)
$$
where $\phi(F)$ denotes half the sum of the framings on $L$ induced by
$F$.  (Note that these framings are all even since $L$ is proper.)

\examples{(1)\qua The Borromean rings have Brown invariant 4.  This can be
seen using the bounded checkerboard surface $F$ in the minimal diagram for
the link shown in Figure 2a.  The enhanced homology is $P_1\oplus 2A_0$, and
the induced framings are all $-2$, so the Brown invariant is
$1-{\frac{1}{2}}(-6) = 4$.

\smallskip

(2)\qua The $k$--twisted Bing double of any knot (with $k$ full twists in the
parallel strands) has Brown invariant $4k$.  To see this use the
obvious banded Seifert surface of genus 1 (shown for the double of the
unknot in Figure 4b) which has enhancement $T_{4k}\oplus A_0$ and
induces the 0--framing on both boundary components.}

{\small\vskip .1in
\centerline{\fig{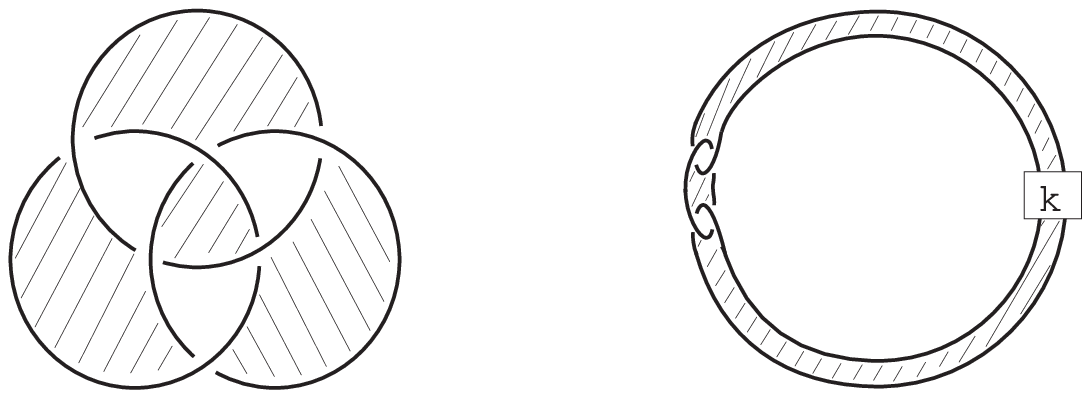}{.8}}
\vskip .1in
\centerline{\quad (a) Borromean rings \kern 60pt (b) $k$-twisted Bing
double}
\figno{4}}

\remarks{(1)\qua If $L$ is given an orientation, then it has an Arf
invariant $\alpha(L) \in \bz_2$ which is related to the Brown invariant of
any (oriented) Seifert surface $F$ for $L$ by the identity $\beta(F) =
4\alpha(L)$.  Adjusting for the framings one obtains the formula
$$\beta(L) = 4\,\alpha(L) + \lk(L)$$
where $\lk(L)$ denotes the sum of all the pairwise linking numbers of $L$.
Note that both terms on the right hand side depend on the orientation of
$L$, while their sum does not.  For example the $(2,4)$-torus link $L$ has
$\beta(L)=6$, while $(\alpha(L),\lk(L)) = (1,2)$ or $(0,-2)$ according to
whether the components are oriented compatibly or not.


(2)\qua (see \cite{kt})  The Brown invariant of $L$ can also be defined using
a surface $F$ in $B^4$ bounded by $L$ for which there exists a Pin$^-$
structure on $B^4 - F$ which does not extend over $F$.  The Pin$^-$
structure descends to a Pin$^-$ structure on $F$ which determines an
element in 2--dimensional Pin$^-$ bordism which, using the Brown
invariant, is isomorphic to $\bz_8$.}

For an immersion $f\co F\looparrowright S^3$ bounded by $L$, a correction
term coming from the singularities of $f$ is needed to compute $\beta(L)$
in terms of $\beta_f$.  This is most easily expressed using the {\it
quarter-twist map}
$$q\co\{\textup{immersed doublebands in $S^3$}\} \to \bz_8,$$
defined analogously to the half-twist map $h$ above:  A {\it doubleband}
is a union of $\times$--bundles over circles; for embedded doublebands
$B$ with core $C$, define $q(B) = \lk(C,\partial B) \pmod8$, where
$\partial B$ is the compatibly oriented $::$--bundle ($S^0\times S^0-$bundle)
and then observe that this is a regular homotopy invariant.   For example, any
double curve $C$ in $F' = f(F)$ has an immersed doubleband neighborhood
$B$, and $q(B)$ (also denoted $q(C)$ by abuse of notation) is odd or even
according to whether $f^{-1}(C)$ consists of one or two curves in $F$; we
say $C$ is {\it orientation-reversing} in the former case, and {\it
orientation-preserving} in the latter.

Now consider the entire singular set of $f$.  It consists of a
collection of double curves which intersect in some number $\tau_f$ of
triple points.  A neighborhood of this singular set is an immersed
doubleband $B$ (generally not connected) and we define $\delta_f =
q(B)$, and (as for the case of embedded surfaces) $\phi_f$ to be half the sum of
the framings on $L$ induced by $F$.

\begin{theorem}
Let $L$ be a proper link bounded by a regularly immersed surface
$f\co F\looparrowright S^3$.  Then the Brown invariant
$$\beta(L) = \beta_f  - \phi_f + 3\delta_f +4\tau_f,$$
and so (by Remark 1 above) the Arf invariant
$$\alpha(L) = \textstyle{\frac{1}{4}}(\beta_f - \phi_f - \lk(L)
  + 3\delta_f +4\tau_f)$$
for any chosen orientation on $L$.
\end{theorem}

\begin{proof}
The strategy is to reduce to the embedded case by local modifications of
$f$.  We first eliminate triple points by {\it Borromean cuts} as
shown in Figure 5.   This calls for the removal of three disks (bounded by
the Borromean rings) and the addition of three tubes to maintain {\sl
interior} singularities (a condition for regularity of the immersion) in the
three sheets near each triple point.  The effect on the boundary is to add
$\tau_f$ copies of the Borromean rings, which changes the left hand side of
the formula by $4\tau_f$ (see Example 1 above).  The terms $\beta_f$ and
$\phi_f$ on the right hand side are unchanged, since the effect of this
modification is to add copies of $T_0\oplus A_0$ to the enhanced homology,
and to give the 0--framing to the Borromean rings on the boundary.
Likewise $\delta_f$ is unchanged, since the double curves have simply
undergone a regular homotopy, so the net change on the right hand side
is also $4\tau_f$.

{\small\vskip .1in
\centerline{\fig{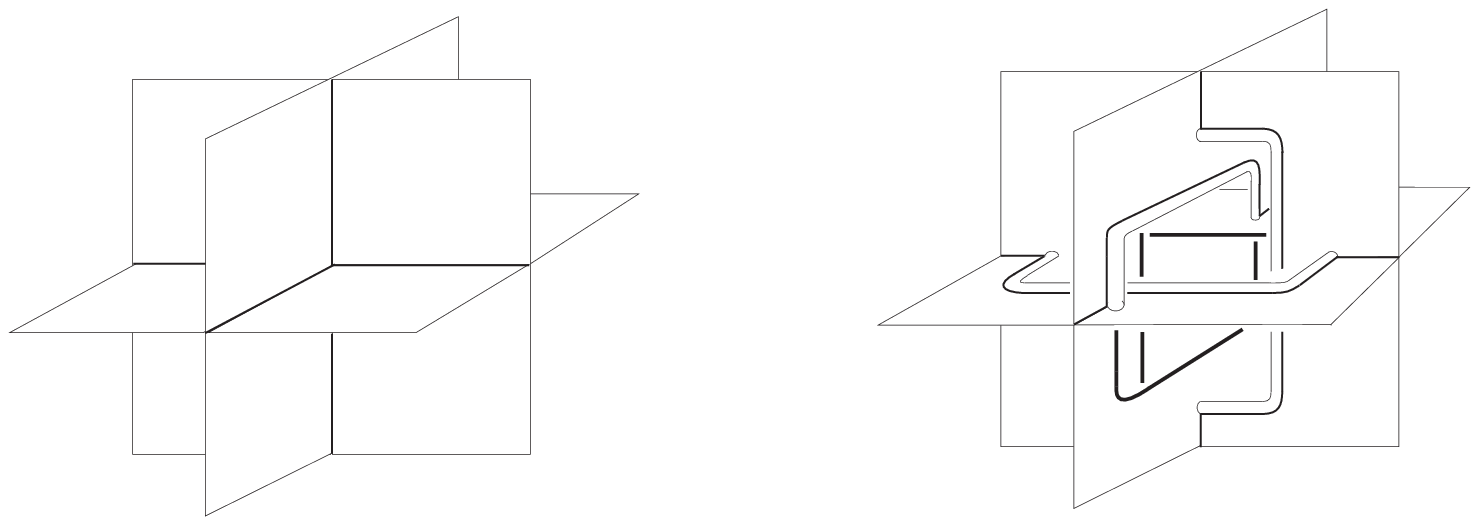}{.8}}
\centerline{\quad (a) triple point \kern 100pt (b) Borromean cut}
\figno{5}}

Now fix a double curve $C$, and set $n=q(C)$, the number of quarter-twists
in the normal $\times$-bundle of $C$.  Note that $C$ is embedded since
there are no triple points.  Let $g\co E\looparrowright S^3$ be the immersion
obtained from $f$ by ``smoothing'' along $C$.  In other words, proceed along
$C$, replacing each fiber in the normal $\times$--bundle of $C$ (see Figure
6a) by two arcs (Figure 6b); if $C$ is orientation reversing (ie $n$
is odd) then one must insert a saddle at some point to allow the fibers to
match up (Figure 6c).  This can be done in two ways, starting with
\break \smooth \ or $)\,\!($, and either will do.  In any case $C$ disappears
from set of double curves, and so $\delta_g = \delta_f-n$.  Since this
construction adds no new boundary components or triple points, it suffices to
show that
$\beta_g = \beta_f+3n$.   There are several cases to consider depending on
the value of $n$.

{\small\vskip .3in
\centerline{\fig{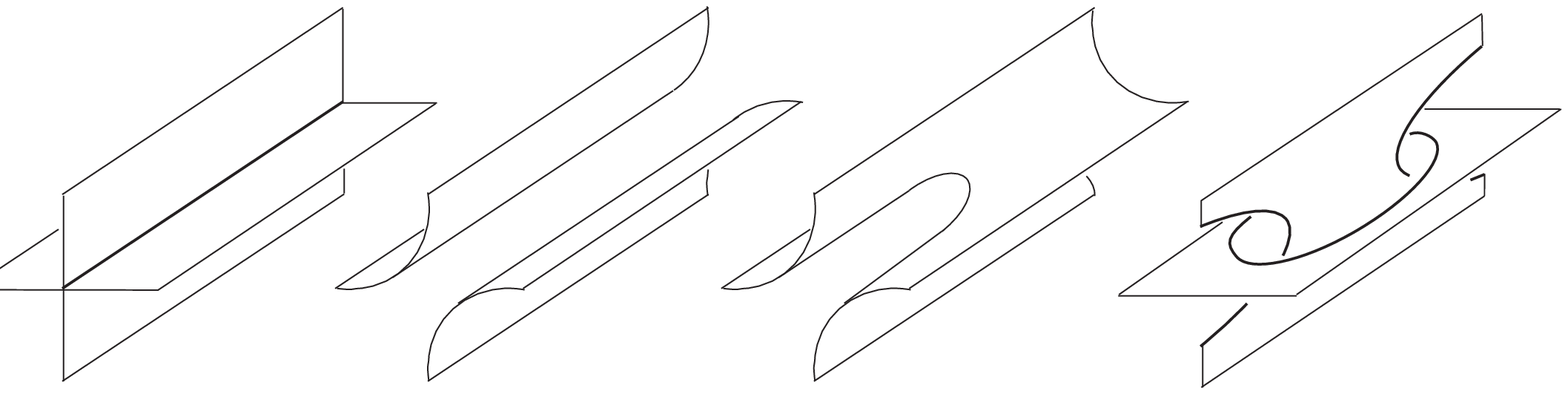}{.60}}
\centerline{(a) double curve \place{-25}{25}{$F$} \kern 6pt (b) even
smoothing \place{-50}{25}{$E$}
\kern 6pt (c) odd smoothing \place{-65}{25}{$E$}
\kern 6pt (d) Bing cut}
\figno{6}}

If $n = 4k\pm1$, for some $k$, then $E$ is obtained from $F$ by removing
a M\"obius band, whose image wraps twice around $C$, and replacing it with
a boundary-connected sum of two M\"obius bands.  The effect on the enhanced
homology is to delete a $P^{n+2}$--summand (note that a regular homotopy of
the doubly wrapped M\"obius band introduces a kink which adds 2 to the
number of half-twists) and to add two $P^{2k\pm1}$--summands.  As above, the
superscript is the value of the enhancement on a generator, and so the
Brown invariant is obtained by reducing mod 4 to $\pm1$.  In other words,
we are deleting a $P_{\mp1}$-summand and adding two $P_{(-1)^k}$--summands.
Thus $\beta_g \equiv \beta_f\pm1\pm2(-1)^k \equiv \beta_f+3n \pmod8$.

If $n = 4k+2$, then $E$ is obtained by removing two M\"obius bands from
$F$ and then identifying the resulting boundary components to form a circle
$C'$.  The effect on the enhanced homology is to delete two $P^{n/2} =
P_{(-1)^k}$--summands, and to add either nothing (if the M\"obius bands lie
in distinct components of $F$) or one $T_0$ or $K_0$--summand\foot{By
definition $K_0=P_1\oplus P_{-1}$.  This corresponds to a Klein bottle $K
\cong P\oplus P$.  The four enhancements $K^{0,\pm1}, K^{2,\pm1}$ of $K$,
with respect to the homology basis $x,y$ with $x\cdot x = 0$ and $x\cdot y
= y\cdot y = 1$, fall into three isomorphism classes $K_0 = \{K^{0,1},
K^{0,-1}\}$, and $K_{\pm2} = K^{2,\pm1}$.} (a punctured torus or Klein
bottle is seen as a neighborhood of the union of $C'$ and a dual circle).
Thus $\beta_g \equiv \beta_f-2(-1)^k \equiv \beta_f-n \equiv \beta_f+3n
\pmod8$.

Finally consider the case $n = 4k$.  Let $C_1,C_2$ denote the two
circles in $f^{-1}(C)$.  An argument similar to the one above can be given
by analyzing several cases depending on the homology classes of $C_1$,
$C_2$ and $C_1\cup C_2$, but there is a simpler argument using a different
modification of $f$ near $C$ which we call a {\it Bing cut}.  It is
obtained by removing two discs from $F$, thereby introducing a 0--framed
$k$--twisted Bing double of $C$ on the boundary (see Example 2 above) as
shown in Figure 6d.  This adds $4k$ to $\beta(L)$ since $C$ links $L$
evenly; indeed $\lk(C,L) \equiv \lk(C,L\cup C^+\cup C^-) \equiv 0\pmod2$,
where $C^\pm$ are pushoffs of $C$ in the image of a neighborhood of $C_1$,
since $L\cup C^+\cup C^-$ bounds in $S^3-C$ across $f(F-C_2)$.  The terms
on the right hand side remain unchanged except for $3\delta_f$ which also
changes by $12k \equiv 4k \pmod8$.  The proof is completed by induction on
the number of double curves.
\end{proof}

A formula for the Brown invariant of a totally proper link $C$ can be
deduced by applying this theorem to a family of connected surfaces (bounding
the components $C_i$ of $C$) which meet only in their interiors.  In 
this case, $\beta_f-\phi_f = \sum \beta(C_i)$, $\delta_f = \sum 
\lambda(C_i,C_j)$ and $\tau_f = \sum \tau(C_i,C_j,C_k)$, where 
$\lambda$ and $\tau$ are the Sato--Levine and Milnor invariants 
defined in section 1.  Noting that $\lambda$ is even valued, we have
$$\beta(C) = \sum_i\beta(C_i) - \sum_{i<j} \lambda(C_i,C_j) +
  4\sum_{i<j<k} \tau(C_i, C_j, C_k) \pmod 8.$$
The formula 
$$\alpha(C) = \sum_i\alpha(C_i) + \sum_{i<j}
  {\textstyle{\frac{1}{4}}}(\lambda(C_i,C_j) + \lk(C_i,C_j)) +
  \sum_{i<j<k} \tau(C_i, C_j, C_k) \pmod 2$$
for the Arf invariant of $C$ with any chosen orientation (Theorem 
1.1) follows, since $\beta(C) = 4\alpha(C) + \lk(C)$.  The dependence 
on the orientation is captured by $\lk(C)$, the sum of the pairwise 
linking numbers of the components of $C$.

\Addresses

\end{document}